\documentclass[11pt]{article}
\usepackage{dsfont}

\usepackage{amsmath}
\usepackage{float}
\usepackage{hyperref}
\usepackage{mathrsfs}
\usepackage{amssymb}
\usepackage{caption}
\usepackage{subcaption}
\usepackage{tikz}
\usepackage{verbatim}
\usepackage{color}
\tikzstyle{every node}=[circle,inner sep=1pt]
\usepackage{graphicx}
\usepackage{times}
\usepackage{cite}

\newtheorem{theorem}{Theorem}[section]

\newtheorem{proposition}[theorem]{Proposition}

\usetikzlibrary{trees}

\tikzstyle{every node}=[circle,inner sep=2pt]
\tikzstyle{tn}=[shape=circle, draw, color=black!70]
\tikzstyle{tn1}=[shape=circle,fill=white!60, draw, color=black!70, inner sep=1.5pt]
\tikzstyle{tn2}=[shape=circle,font=\fontsize{2}{0}\selectfont, draw, color=black!70, inner sep=0.6pt]
\tikzstyle{marke}=[shape=circle,minimum size=0.2cm, draw,blue]

\makeatletter \@addtoreset{equation}{section} \makeatother
\makeindex 
\def\qed{\hfill \rule{4pt}{7pt}}

\begin{document}

\begin{center}
{\bf The Gessel Correspondence and the Partial $\gamma$-Positivity of

the Eulerian Polynomials
  on Multiset Stirling Permutations}

\vskip 6mm

William Y.C. Chen$^1$, Amy M. Fu$^2$
and Sherry H.F. Yan$^3$

\vskip 3mm

$^{1}$Center for Applied Mathematics\\
Tianjin University\\
Tianjin 300072, P.R. China

\vskip 3mm

$^{2}$School of Mathematics\\
Shanghai University of Finance and Economics\\
Shanghai 200433, P.R. China

\vskip 3mm

$^3$Department of Mathematics\\
Zhejiang Normal University\\
Jinhua, Zhejiang 321004, P.R. China

Emails: { chenyc@tju.edu.cn, fu.mei@mail.shufe.edu.cn, hfy@zjnu.cn }

\end{center}

\begin{center}
    {\bf Abstract}
\end{center}

Pondering upon the grammatical labeling
of 0-1-2 increasing plane trees, we come to the
realization that
the grammatical labels play a role
as records of chopped off leaves of
 the original increasing binary trees. While such an
understanding is purely psychological, it does
give rise to an efficient apparatus to tackle
the partial $\gamma$-positivity
of the Eulearian polynomials on
multiset Stirling permutations, as long as
we bear in mind the combinatorial meanings of
the labels $x$ and $y$  in the Gessel  representation of
a $k$-Stirling permutation by means of an
increasing $(k+1)$-ary tree. More precisely, we introduce a Foata-Strehl
action on the Gessel trees resulting in an interpretation
of the partial $\gamma$-coefficients of the aforementioned
Eulerian polynomials,  different from
the ones found by Lin-Ma-Zhang and
Yan-Huang-Yang. In particular, our strategy
can be adapted to deal with
the partial  $\gamma$-coefficients of the
second order Eulerian polynomials, which in turn can be
readily converted to
the combinatorial formulation due to Ma-Ma-Yeh in connection
with certain statistics of Stirling permutations.



\section{Introduction}

This work is concerned with the partial
$\gamma$-coefficients of the
Eulerian polynomials on multiset
Stirling permutations,
which are also called the Stirling polynomials.

For $n\geq 1$, let  ${S}_n$ denote the set of
permutations of $[n]=\{1, 2, \ldots, n\}$.
For a permutation $\sigma=\sigma_1\sigma_2\cdots \sigma_n$
in $S_n$, we adopt the convention that
a zero is patched both at the beginning and
at the end of $\sigma$, that is,
$\sigma_0=\sigma_{n+1}=0$.
An ascent of $\sigma$ is
defined to be an index $i$ $(1\leq i \leq n)$
such that $\sigma_{i-1}<\sigma_i$,
whereas a descent is
defined to be an index $i$ $(1\leq i \leq n)$
such that $\sigma_i > \sigma_{i+1}$.
The numbers of ascents and descents of
$\sigma$ are denoted by ${\rm asc}(\sigma)$
and ${\rm des}(\sigma)$, respectively.
The
bivariate Eulerian polynomials $A_n(x, y)$ are defined by
\begin{equation}\label{A_n}
A_n(x,y)=\sum_{\sigma\in {S}_n} x^{{\rm asc}(\sigma)}y^{{\rm des}(\sigma)}.
\end{equation}
Setting $y=1$, $A_n(x,y)$ takes the form
of the usual Eulerian polynomials, or
the descent polynomials of $S_n$.

One of the most  remarkable facts about the Eulerian polynomials is the
$\gamma$-positivity discovered by
Foata and Sch\"uzenberger \cite{FS-70},
which has been extensively studied ever since,
see, for example, \cite{A-2018, Branden-2008, Chow-2008, Elizalde-2021}.

We shall choose to work with the
bivariate version $A_n(x,y)$.
The following expression of
$A_n(x, y)$ is called $\gamma$-expansion:
\begin{equation}\label{gammaA}
A_n(x,y)= \sum_{k=1}^{\lfloor (n+1)/2\rfloor} \gamma_{n,k} (xy)^k (x+y)^{n+1-2k}.
\end{equation}
The coefficients $\gamma_{n,k}$ are
called the $\gamma$-coefficients.
Foata and Sch\"uzenberger
discovered a combinatorial interpretation of the $\gamma$-coefficients
implying the positivity. More precisely, it has been shown that
  \begin{equation}
   \gamma_{n,k}=|\{\sigma\in   {S}_n\mid {\rm des}(\sigma)=k,
 {\rm ddes}(\sigma)=0\}|,
 \end{equation}
 where ${\rm ddes}(\sigma)$ means the number of double descents of $\sigma$,
 that is, the number of indices $i$
 such that $\sigma_{i-1}> \sigma_i >\sigma_{i+1}$.

As a notable extension of the Eulerian polynomials,
 Gessel and Stanley
\cite{Gessel-Stanley-1978} introduced the
notion of Stirling permutations whose descent
polynomials have been called the second order
Eulerial polynomials.

For $n\geq 1$, let $[n]_2$ denote
the multiset $\{1^2, 2^2, \ldots, n^2\}$,
where $i^2$ signifies two occurrences of $i$. A permutation $\sigma$
on $[n]_2$ is said to be a Stirling
permutation if for any $i$,
the elements between the two occurrences
of $i$ in $\sigma$, if any, are greater than $i$. For $n\geq 1$, the set of
  Stirling permutations of $[n]_2$ is usually denoted by $Q_n$. As before,
  we assume that a Stirling permutation
  is patched a zero both at the beginning
  and at the end. The statistics ${\rm asc}$ and ${\rm des}$ can be analogously defined for Stirling permutations.

    The number of Stirling permutations in $Q_n$ with $k+1$ descents,
    often denoted by $C(n,k)$, is called the second order Eulerian number.
    The generating function $C_n(x)$ of $C(n,k)$ has been referred to as
    the second order Eulerian polynomial.
        For Stirling permutations,
    one more statistic naturally comes on the scene, that is, the number of plateaux. It appears that the notion
    of a plateau was first introduced by
    Dumont \cite{Dumont-1980} under the name of
    a repitition. Let $\sigma=\sigma_1\sigma_2\cdots \sigma_{2n} \in Q_n$. An index $i$
    $(1\leq i \leq 2n)$ is called a plateau
    if $\sigma_{i}=\sigma_{i+1}$. The
    number of plateaux of $\sigma$ is
    denoted by ${\rm plat}(\sigma)$.

   B\'ona \cite{Bona-2008}  proved that   the three
statistics ${\rm asc}$, ${\rm plat}$ and ${\rm des}$
are equidistributed over $Q_n$.
 Janson \cite{Janson-2008}   constructed an urn model to prove the symmetry of the joint distribution of the three statistics.

In fact, Dumont \cite{Dumont-1980}
defined the trivariate second order Eulerian polynomials
\begin{equation}\label{C-n}
 C_{n}(x,y,z)=\sum_{\sigma\in {Q}_{n}}
 x^{{\rm asc}(\sigma)}y^{{\rm des}(\sigma)}z^{{\rm plat}(\sigma)},
  \end{equation}
 which can be regarded as an extension
 to the second order Eulerian polynomials
 of Gessel-Stanley and the bivariate
 Eulerian polynomials. Apparently, when
 a Stirling permutation $\sigma \in Q_n$
 has $n$ plateaux, it can be considered
 as a permutation on $[n]$ with each element $i$ replaced by $ii$.
   It was noticed by Dumont that
  $C_n(x,y,z)$
  are symmetric in $x,y,z$.

 The question of $\gamma$-positivity
 for $C_n(x,y,z)$ has been studied by
 Ma-Ma-Yeh \cite{Ma-Ma-Yeh-2019}. Write
\begin{equation}\label{c-n-gamma-1}
    C_n(x,y,z) = \sum_{i=1}^n z^i \sum_{j=0}^{\lfloor (2n+1-i)/2\rfloor}
      \gamma_{n, i, j} (xy)^j (x+y)^{2n+1-i-2j},
\end{equation}
which is called the partial $\gamma$-expansion
of $C_n(x,y,z)$. The coefficients $\gamma_{n,i,j}$
are called the partial $\gamma$-coefficients.
Making use of a context-free grammar argument,
Ma-Ma-Yeh showed that $C_n(x,y,z)$ are partial $\gamma$-positive
in the sense that the coefficients $\gamma_{n,i,j}$ are nonnegative.
Moreover, they obtained a combinatorial interpretation
of $\gamma_{n,i,j}$ resorting to certain statistics on
Stirling permutations.

The structure of a Stirling permutation can be
further extended to a multiset. Throughout this paper,
we assume that $n\geq 1$. Unless specified otherwise, we always assume that
\begin{equation}\label{m}
M = \{ 1^{k_1}, 2^{k_2}, \ldots, n^{k_n} \},
\end{equation}
where $k_i\geq 1$ for all $i$ and $i^{k_i}$
stands for $k_i$ occurrences of $i$.
Moreover, we always designate $K$ to denote $k_1+k_2+\cdots + k_n$.

   A permutation $\sigma=\sigma_1\sigma_2\ldots \sigma_K$ of  $M$
  is said to be  a {  Stirling } permutation if  $\sigma_i=\sigma_j$ with $i<j$, then $\sigma_k \geq \sigma_i$ for any $i<k<j$.
  The set of Stirling permutations of $M$
  will be denoted by $Q_M$.
  For $M=[n]_k=\{ 1^k, 2^k, \ldots, n^k\}$,
  a Stirling permutation on $M$ is called
  a $k$-Stirling permutation.

The statistics   ${\rm asc} $, ${\rm des} $ and  ${\rm plat}$ for Stirling permutations
in $Q_n$ can be literally carried over to
$Q_M$. Then
the Eulerian
polynomials
on  Stirling permutations of a multiset $M$ can be defined by
 \begin{equation}\label{Q_n}
 C_{M}(x,y,z)=\sum_{\sigma\in {Q}_{M}}x^{{\rm asc}(\sigma)}y^{{\rm des}(\sigma)}z^{{\rm plat}(\sigma)},
  \end{equation}
 see also \cite{Lin-Ma-Zhang-2021}. In particular, the descent polynomial
 over $Q_M$ is often denoted by
 \begin{equation}
     Q_M(x)= \sum_{\sigma\in {Q}_{M}}x^{{\rm des}(\sigma)},
 \end{equation}
 see also \cite{Yan-Zhu-2022}. As shown by Frenti \cite{Brenti-1989},
 $Q_M(x)$ has only real roots for any multiset $M$.

 While $C_M(x,y,z)$ are no longer symmetric
 in general, they are symmetric in $x$ and $y$. This means that
   $C_{M}(x,y,z)$ can be expressed as
 \begin{equation}\label{gammaQ}
C_{M}(x,y,z)=\sum_{i=0}^{K-n}z^i \sum_{j=1}^{\lfloor (K+1-i)/2 \rfloor} \gamma_{M, i,j}(xy)^j(x+y)^{K+1-i-2j}.
 \end{equation}
The above relation (\ref{gammaQ}) is
called the { partial $\gamma$-expansion} of $C_{M}(x,y,z)$. The coefficients $\gamma_{M,i,j}$ are called the { partial $\gamma$-coefficients} of  $C_{M}(x,y,z)$. The nonnegativity of the  coefficients  $\gamma_{M, i,j}$ is referred to as the {partial $\gamma$-positivity}.

The  partial $\gamma$-positivity
for several multivariate ploynomials associated various classes of permutations
has recently been studied in
Ma-Ma-Yeh \cite{Ma-Ma-Yeh-2019},
Lin-Ma-Zhang \cite{Lin-Ma-Zhang-2021} and  Yan-Huang-Yang \cite{Yan-Huang-Yang-2021}.

The objective  of this paper
is to present a combinatorial treatment   of
the  partial $\gamma$-coefficients of
$C_M(x,y,z)$.
In particular, our strategy
can be adapted to deal with
the  partial $\gamma$-coefficients of the
second order Eulerian polynomials, which in turn can be
readily converted to
the combinatorial formulation obtained
by Ma-Ma-Yeh \cite{Ma-Ma-Yeh-2019}.
  Our  combinatorial interpretation
of $\gamma_{M,i,j}$  is built on the Gessel trees
which are increasing plane trees in which the  internal vertices are
represented by distinct numbers,
whereas in
the combinatorial framework of Yan, Huang and Yang
\cite{Yan-Huang-Yang-2021}, two vertices
are allowed to be represented by the
same number.

The underlying combinatorial structure
that is concerned with in this work is that of
a Gessel tree for a Stirling permutation of a multiset.
Our main result (Theorem \ref{mainthT}) is a combinatorial interpretation
in light of canonical Gessel trees. A closely related approach is
to utilize context-free grammars, which leads to the notion of pruned Gessel
trees. {A careful study of the Gessel correspondence
reveals the  properties  needed to
turn Theorem \ref{mainthT}
into an equivalent
statement  on multiset
Stirling permutations (Theorem \ref{mainthQ}).}
Specializing to the set  $Q_n$ of
Stirling permutations, we are led to the combinatorial interpretation
of the partial $\gamma$-coefficients (Theorem \ref{thma})
due to Ma-Ma-Yeh \cite{Ma-Ma-Yeh-2019}.

The rest of this paper is organized as follows. In Section 2,
we give an overview of the Gessel correspondence
between increasing trees and multiset  Stirling permutations. We also address a refined property of the Gessel correspondence. Section 3 is devoted to
a classification of Gessel trees based
on the structure of canonical Gessel trees. An operation in the spirit of the
Foata-Strehl action \cite{Foata-Strehl-1974}
is introduced to serve the purpose. As a main result of this work, we present a combinatorial
interpretation of the partial $\gamma$-coefficients of $C_M(x,y,z)$. In Section 4, we present a context-free grammar approach. In fact, such a consideration has spurred the notion of
pruned Gessel trees and has provided a
motivation for the Foata-Strehl action in Section 3. The aim of Section 5 is to present a combinatorial explanation
of the
partial $\gamma$-coefficients in terms
of multiset Stirling permutations. In Section 6, we explain how to get to the
result of Ma-Ma-Yeh. It turns out that the
symmetry property of $C_n(x,y,z)$ is
required to fulfill the task.

\section{The Gessel correspondence}

The Gessel correspondence \cite{Gessel-2020} is an
extension of the classical
bijection between permutations and
increasing binary trees to Stirling permutations on a multiset, see also \cite{Chow-2008,
JKP-2011, Elizalde-2021}.

It is the aim of this paper to
utilize this correspondence to
study the $\gamma$-positivity and
the partial $\gamma$-positivity of the Eulerian polynomials on
multiset Stirling permutations.
Gessel   \cite{Gessel-2020}  established a bijection
$\phi$ between $k$-Stirling permutations and $(k+1)$-ary
increasing trees in which only the
internal vertices are labeled and represented by
solid dots, whereas the external
vertices (leaves) are not labeled and represented  by   circles.

 For a Stirling permutation $\sigma$ on $M$. If $\sigma=\emptyset$, let $\phi(\sigma)$ be the tree with only one (unlabeled) vertex.  If  $\sigma\neq \emptyset$, let $i$ be the smallest element of $\sigma$. Then $\sigma$ can be uniquely decomposed as $w_0 \,i \, w_1 \, i \, \cdots \, w_{k-1} \,
 i \, w_k$, where $w_j$ is either empty or a Stirling permutation for all $0 \leq j\leq k$.
Set  $i$ to be the root of $\phi(T)$,
and set  $\phi(T_j)$ to be the $j$-th subtree (counting from left to right) of $i$. This procedure yields a recursive construction
 of $\phi(\sigma)$.

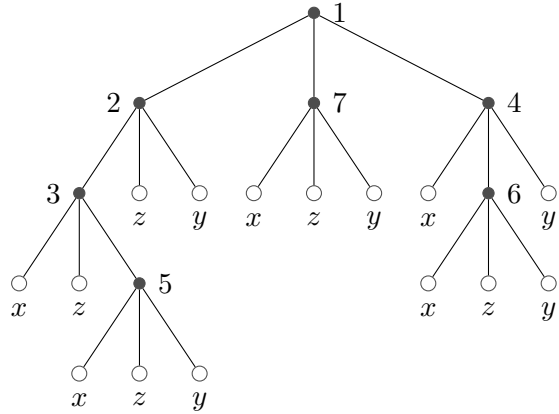
\begin{figure}[!ht]
\begin{center}
\begin{tikzpicture}[scale=0.8]
\node [tn1, label=0:{$1$}]{}[grow=down]
	[sibling distance=29mm,level distance=15mm]
        child {node [tn1,label=180:{$2$}](five){}
            [sibling distance=10mm,level distance=15mm]
             child{node [tn1,label=180:{$3$}]{}
            [sibling distance=10mm,level distance=15mm]
            child{node [tn,label=-90:{$x$}]{}}
            child{node [tn,label=-90:{$z$}](four){}}
         child{node [tn1,label=0:{$5$}](three){}
             [sibling distance=10mm,level distance=15mm]
                child{node [tn,label=-90:{$x$}](five1){}}
            child{node [tn,label=-90:{$z$}](five2){}}
            child{node [tn,label=-90:{$y$}]{}
              edge from parent
                    node[right=2pt ]{}
              }
              edge from parent
                    node[right=2pt ]{}
                    }
             edge from parent
                    node[right=2pt ]{}
                    }
            child{node [tn,label=-90:{$z$}](four){}}
            child{node [tn,label=-90:{$y$}](three){}
            edge from parent
                node[left=2pt ]{}
                }
            edge from parent
                node[left=2pt ]{}
            }
        child {node [tn1,label=0:{$7$}]{}
        [sibling distance=10mm,level distance=15mm]
             child{node [tn,label=-90:{$x$}]{}}
             child{node [tn,label=-90:{$z$}]{}}
             child{node [tn,label=-90:{$y$}]{}
                edge from parent
                node[left=1pt ]{}
                }
            edge from parent
                node[left=1pt ]{}
                }
        child {node [tn1,label=0:{$4$}]{}
           [sibling distance=10mm,level distance=15mm]
                child{node [tn,label=-90:{$x$}](three){}}
            child{node [tn1,label=0:{$6$}](four){}
                [sibling distance=10mm,level distance=15mm]
                child{node [tn,label=-90:{$x$}]{}}
            child{node [tn,label=-90:{$z$}]{}}
            child{node [tn,label=-90:{$y$}]{}
              edge from parent
                    node[right=2pt ]{}
              }
            edge from parent
                    node[right=2pt ]{}
            }
            child{node [tn,label=-90:{$y$}]{}
              edge from parent
                    node[right=2pt ]{}
              }
            edge from parent
                node[left=2pt ]{}
                }
        ;

\end{tikzpicture}

\end{center}
\caption{A  ternary
 increasing tree.}\label{f1}
\end{figure}

For example,
let $\sigma= 33552217714664  \in {Q}_{7}$, the corresponding  ternary
 increasing tree is illustrated in Figure \ref{f1}.

 Recall that for the classical representation of
 permutations by  increasing
binary trees, every vertex (regardless
of an internal vertex or a leaf) is labeled,
whereas in the Gessel representation, only the
internal vertices are labeled, and external
vertices (leaves) are not labeled, where a leaf
is drawn as a circle.

In the Gessel  representation, a leaf
is called an {$x$-leaf} if it is the first
child, and is called a {$y$-leaf} if it is
the last child; otherwise, it is called a
$z$-leaf. Such leaves are crucial
for the study of the $\gamma$-positivity. As will be seen,
such leaves are the natural ingredients of a
Foata-Strehl action.

Furthermore, Janson-Kuba-Panholzer \cite{JKP-2011} introduced the notion of a
$j$-plateau of a $k$-Stirling permutaiton.
In this event, an index $i$ is said to be a
{$j$-plateau} of $\sigma$
if $\sigma_i=\sigma_{i+1}=r$ and $\sigma_{i}$
is the $j$-th occurrence of $r$ in $\sigma$.
If a leaf $v$ of $T$ is the $j$-th child for some $2\leq j\leq k$, then $v$ is called a { $z_j$-leaf}. Janson-Kuba-Panholzer \cite{JKP-2011} showed that a  $z_j$-leaf in
a $(k+1)$-ary increasing tree
corresponds to a $j$-plateau of  $\sigma$.

\begin{figure}[!ht]
\begin{center}
\begin{tikzpicture}[scale=0.75]
\node [tn1,label=0:{$1$}]{}[grow=down]
	[sibling distance=18mm,level distance=15mm]
        child {node [tn1,label=180:{$2$}](five){}
            [sibling distance=19mm,level distance=15mm]
            child{node [tn1,label=180:{$3$}]{}
            [sibling distance=10mm,level distance=15mm]
              child{node [tn1,label=180:{$5$}](three){}
             [sibling distance=10mm,level distance=15mm]
                child{node [tn,label=-90:{$x$}](five1){}}
            child{node [tn,label=-90:{$z$}](five2){}}
            child{node [tn,label=-90:{$y$}]{}
              edge from parent
                    node[right=2pt ]{}
              }
              edge from parent
                    node[right=2pt ]{}
                    }
            child{node [tn,label=-90:{$z$}](four){}}
            child{node [tn,label=-90:{$y$}]{}
            }
             edge from parent
                    node[right=2pt ]{}
                    }
                     child{node [tn,label=-90:{$y$}](three){}
            edge from parent
                node[left=2pt ]{}}
            }
        child {node [tn,label=-90:{$z$}]{}
            edge from parent
                node[left=1pt ]{}
                }
        child {node [tn1,label=0:{$4$}]{}
           [sibling distance=12mm,level distance=15mm]
                child{node [tn,label=-90:{$x$}](three){}}
            child{node [tn1,label=0:{$6$}](four){}
                [sibling distance=8mm,level distance=15mm]
                child{node [tn,label=-90:{$x$}]{}}
            child{node [tn,label=-90:{$z$}]{}}
            child{node [tn,label=-90:{$z$}]{}}
            child{node [tn1,label=0:{$7$}]{}
             [sibling distance=22mm,level distance=15mm]
                child{node [tn,label=-90:{$x$}](five1){}}
            child{node [tn,label=-90:{$y$}]{}
              edge from parent
                    node[right=2pt ]{}
              }
              }
            edge from parent
                    node[right=2pt ]{}
            }
            child{node [tn,label=-90:{$y$}]{}
              edge from parent
                    node[right=2pt ]{}
              }
            edge from parent
                node[left=2pt ]{}
                } ;
\end{tikzpicture}
\end{center}
\caption{A Gessel tree on $M =\{1^2, 2,3^2,4^2,5^2, 6^3, 7\}$. }
\label{f2}
\end{figure}
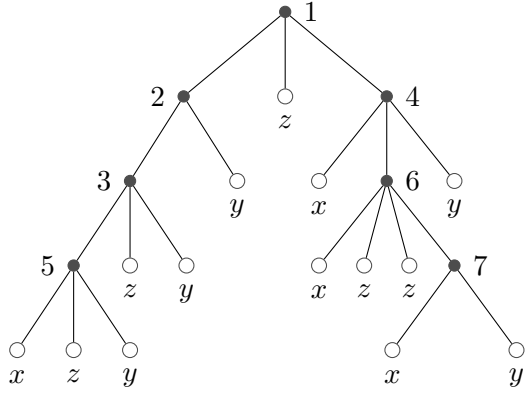

By a Gessel tree on $M$ we mean a
plane tree with internal vertices $1, 2, \ldots, n$
along with unlabeled external vertices such that the internal
vertex $i$ has exactly $k_i+1$
children and the internal vertices form an increasing tree.
For instance,   Figure \ref{f2} exhibits a Gessel
tree on $M =\{1^2, 2,3^2,4^2,5^2, 6^3, 7\}$.

 Denote by ${\rm xleaf}(T)$, ${\rm yleaf}(T)$ and ${\rm zleaf}(T)$ the numbers of $x$-leaves, $y$-leaves and $z$-leaves, respectively.
The following property can be easily deduced from the
recursive construction of the Gessel correspondence.

\begin{proposition}\label{lem1}
Let $n$ and $M$ be given as before.
The Gessel map $\phi$ establishes a one-to-one
correspondence between Stirling permutations  of $M$ and  Gessel
trees on $M$. Moreover,  let
$\sigma\in {Q}_{M}$  and  $T=\phi(\sigma)$. Then
\begin{equation}\label{sigma-T}
({\rm asc}(\sigma), {\rm des}(\sigma), {\rm plat}(\sigma))=({\rm xleaf}(T), {\rm yleaf}(T), {\rm zleaf}(T)).
\end{equation}
\end{proposition}

Instead of reproducing a proof of the above property,
we discuss a  refined description of the Gessel correspondence that
will be needed later in this paper. For this purpose, we shall
introduce the notion of the Gessel decomposition of a
Stirling permutation on $M$.

Let $\sigma$ be a Stirling permutation on $M$.
For any $1\leq i \leq n$, the $i$-segment of $\sigma$, denoted by
$S_i(\sigma)$, is defined to be
the unique sequence $\sigma_r \sigma_{r+1} \cdots \sigma_s$
containing the element $i$, where
$1\leq r \leq s \leq K$, such that
$\sigma_{r-1} < \sigma_{r}$ and $\sigma_{s} > \sigma_{s+1}$
with the convention $\sigma_0=\sigma_{K+1}=0$.
It is clear from the definition of
a Stirling permutation of
$M$ that the $i$-segment of $\sigma$ is well-defined. In fact, one sees that
$S_i(\sigma)$ contains all the
occurrences of $i$ in $\sigma$.

For example, for the Stirling permutation $\sigma=55 33 2 11 4 666 7 4$ of $M =\{1^2, 2,3^2,4^2,5^2, 6^3, 7\}$  corresponding to the Gessel tree in Figure 2, we have
\begin{eqnarray*}
S_1(\sigma)&= & 55 33 2 11 4 666 7 4, \\[1pt]
S_2(\sigma) & = &  55332, \\[1pt]
S_3(\sigma) & = & 5533, \\[1pt]
S_4 (\sigma)& = & 4 666 7 4 , \\[1pt]
S_5(\sigma) & = & 55, \\[1pt]
S_6 (\sigma)&= & 6667 , \\[1pt]
S_7(\sigma) & = & 7.
\end{eqnarray*}

The idea of the Gessel correspondence can be perceived
as a decomposition of an $i$-segment of a Stirling permutation. Let $\sigma$ be
a Stirling permutation on a multiset $M$. The Gessel
decomposition of the $i$-segament $S_i(\sigma)$ is
defined to be a decomposition
\begin{equation}\label{GD} S_i(\sigma)=w_0\,i\, w_1 \, i \, w_2 \, i \, \cdots \, i \, w_{k_i},
\end{equation}
where for each $0\leq t \leq k_i$, $w_t$ is either empty or
a $j$-segment for some $j$.

We now come to the following refined property
of the Gessel correspondence.

\begin{proposition}
\label{RGC}
Let $n$ and $M$ be given as before.
Let $\sigma$ be a Stirling permutation on $M$ and let
$T$ be the corresponding Gessel tree.  Assume that
$\sigma_p$ is the first occurrence of $i$ in $\sigma$ and
$\sigma_q$ is the last occurrence of $i$ in $\sigma$.
Then the vertex $i$ has an $x$-leaf in $T$ if and only if
$p$ is an ascent and $i$ has a
$y$-leaf if and only if $q$
is a descent of $\sigma$.
\end{proposition}

\noindent{\em Proof.}
 First, assume
that $i$ has an $x$-leaf in $T$.
By the Gessel correspondence,
$w_0$ in the Gessel decomposition of $\sigma$ as
given in (\ref{GD}) is empty.
Since an $i$-segment of $\sigma$ is surrounded by two elements smaller than $i$, $\sigma_p$ is
directly preceded by a smaller
element. That means that $p$ is
an ascent of $\sigma$. Conversely, if $p$ is an ascent of $\sigma$, then $w_0$ must be
empty, and so $i$ has an $x$-leaf in $T$. The same reasoning applies to a descent
of $\sigma$ involving the last
occurrence of $i$ and a $y$-leaf of $i$ in $T$, and hence the proof is complete. \qed

 \section{A Foata-Strehl action on
 Gessel trees}

To give a combinatorial interpretation of the $\gamma$-coefficients,
we shall define an  action  on a Gessel tree, which plays the same role
as the original Foata-Strehl group
action for the evaluation of the $\gamma$-coefficients. Such an action
is often called a modified Foata-Strehl action.

Let $T$ be a Gessel tree
on $M$.
We say that an $x$-leaf   in $T$ is balanced if its parent has a $y$-leaf; otherwise, we say that the $x$-leaf is unbalanced. Similarly, a $y$-leaf is said to be
balanced if its parent has an $x$-leaf; otherwise, it is said to be unbalanced.

\begin{figure}[!ht]
\begin{center}
\begin{tikzpicture}[scale=0.75]
\node [tn1,label=0:{$1$}]{}[grow=down]
	[sibling distance=22mm,level distance=15mm]
        child {node [tn1,label=180:{$2$}](five){}
            [sibling distance=21mm,level distance=15mm]
                     child{node [tn,label=-90:{$x$}](three){}}
            child{node [tn1,label=0:{$3$}]{}
            [sibling distance=10mm,level distance=15mm]
              child{node [tn1,label=180:{$5$}](three){}
             [sibling distance=10mm,level distance=15mm]
                child{node [tn,label=-90:{$x$}](five1){}}
            child{node [tn,label=-90:{$z$}](five2){}}
            child{node [tn,label=-90:{$y$}]{}
              edge from parent
                    node[right=2pt ]{}
              }
              edge from parent
                    node[right=2pt ]{}
                    }
            child{node [tn,label=-90:{$z$}](four){}}
            child{node [tn,label=-90:{$y$}]{}
            }
             edge from parent
                    node[right=2pt ]{}
                    }
            }
        child {node [tn,label=-90:{$z$}]{}
            edge from parent
                node[left=1pt ]{}
                }
        child {node [tn1,label=0:{$4$}]{}
           [sibling distance=12mm,level distance=15mm]
                child{node [tn,label=-90:{$x$}](three){}}
            child{node [tn1,label=0:{$6$}](four){}
                [sibling distance=8mm,level distance=15mm]
                child{node [tn,label=-90:{$x$}]{}}
            child{node [tn,label=-90:{$z$}]{}}
            child{node [tn,label=-90:{$z$}]{}}
            child{node [tn1,label=0:{$7$}]{}
             [sibling distance=22mm,level distance=15mm]
                child{node [tn,label=-90:{$x$}](five1){}}
            child{node [tn,label=-90:{$y$}]{}
              edge from parent
                    node[right=2pt ]{}
              }
              }
            edge from parent
                    node[right=2pt ]{}
            }
            child{node [tn,label=-90:{$y$}]{}
              edge from parent
                    node[right=2pt ]{}
              }
            edge from parent
                node[left=2pt ]{}
                } ;
\end{tikzpicture}
\end{center}
\caption{ Action of $\psi_2$
on the Gessel tree in Figure \ref{f2}.}
\label{f3}
\end{figure}
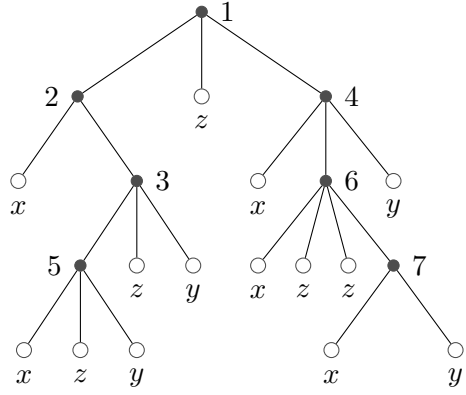

For $1\leq i \leq n$, the Foata-Strehl
action $\psi_i$ is defined as follows. It does nothing to $T$  unless the internal vertex $i$ possesses an unbalanced $y$-leaf. In case
the vertex $i$ has unbalanced $y$-leaf,
then
$\psi_i(T)$ is defined to be the Gessel tree obtained from $T$ by interchanging  the first child
(along with its the subtree) and the last child of vertex $i$, and keeping  the order of  the other children unchanged. As a result, the unbalanced $y$-leaf becomes an unbalanced $x$-leaf. Figure \ref{f3} depicts the action of $\psi_2$ on the Gessel tree in
Figure \ref{f2}.

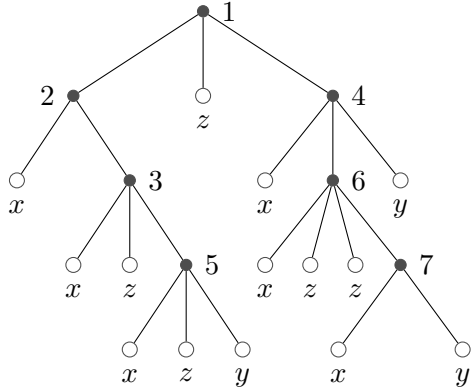
\begin{figure}[h]
\begin{center}
\begin{tikzpicture}[scale=0.75]
\node [tn1,label=0:{$1$}]{}[grow=down]
	[sibling distance=23mm,level distance=15mm]
        child {node [tn1,label=180:{$2$}](five){}
            [sibling distance=20mm,level distance=15mm]
	       child{node [tn,label=-90:{$x$}](three){}}
            child{node [tn1,label=0:{$3$}]{}
            [sibling distance=10mm,level distance=15mm]
            child{node [tn,label=-90:{$x$}](four){}}
            child{node [tn,label=-90:{$z$}]{} }
                child{node [tn1,label=0:{$5$}](three){}
             [sibling distance=10mm,level distance=15mm]
                child{node [tn,label=-90:{$x$}](five1){}}
            child{node [tn,label=-90:{$z$}](five2){}}
            child{node [tn,label=-90:{$y$}]{}
              edge from parent
                    node[right=2pt ]{}
              }
              edge from parent
                    node[right=2pt ]{}
                    }
               edge from parent
                node[left=2pt ]{}
                    }
            edge from parent
                node[left=2pt ]{}
            }
        child {node [tn,label=-90:{$z$}]{}
            edge from parent
                node[left=1pt ]{}
                }
        child {node [tn1,label=0:{$4$}]{}
           [sibling distance=12mm,level distance=15mm]
                child{node [tn,label=-90:{$x$}](three){}}
            child{node [tn1,label=0:{$6$}](four){}
                [sibling distance=8mm,level distance=15mm]
                child{node [tn,label=-90:{$x$}]{}}
            child{node [tn,label=-90:{$z$}]{}}
            child{node [tn,label=-90:{$z$}]{}}
            child{node [tn1,label=0:{$7$}]{}
            [sibling distance=22mm,level distance=15mm]
                child{node [tn,label=-90:{$x$}]{}}
                child{node [tn,label=-90:{$y$}]{}
              edge from parent
                    node[right=2pt ]{}
              }
            edge from parent
                    node[right=2pt ]{}
            }
                        edge from parent
                    node[right=2pt ]{}
            }
            child{node [tn,label=-90:{$y$}]{}
              edge from parent
                    node[right=2pt ]{}
              }
            edge from parent
                node[left=2pt ]{}
                };
\end{tikzpicture}
\end{center}
\caption{A canonical Gessel tree.}
\label{f4}
\end{figure}

Let $G_M$ be the set of Gessel
trees on $M$. By Proposition \ref{lem1},
the Eulerian polynomials $C_M(x,y,z)$
as defined in (\ref{C-n}) can be expressed
in terms of the Gessel trees, namely,
\begin{equation}
    C_M(x,y,z)= \sum_{T\in G_M }x^{ {{\rm xleaf}(T)}}y^{{\rm yleaf}(T)}z^{{\rm zleaf}(T)}.
\end{equation}

The above connection is our starting
point to arrive at  a combinatorial
interpretation of the partial $\gamma$-coefficients of $C_M(x,y,z)$
in the context of Gessel trees.
We need a special class of Gessel trees,
which we call canonical Gessel trees.
To be more specific, we say that
a Gessel tree is canonical if it does not contain any internal vertex having an unbalanced $y$-leaf.

\begin{theorem}\label{mainthT}
Let $n$, $M$ and $K$ be given as before. Then
 \begin{equation}
 \label{gammaT}
 C_{M}(x,y,z)=\sum_{i=0}^{K-n}z^i \sum_{j=1}^{\lfloor (K+1-i)/2 \rfloor} \gamma_{M, i,j}(xy)^j(x+y)^{K+1-i-2j},
 \end{equation}
 where
 $\gamma_{M, i,j}$ is the number of
 canonical Gessel trees on $M$ with $i$   $z$-leaves, $j$ $y$-leaves.
\end{theorem}

\noindent{\em Proof}.
Let $G_{M, i}$ denote  the set of Gessel
trees on $M$  with $i$ $z$-leaves.
 Then Theorem \ref{mainthT} is equivalent to
 \begin{equation}\label{eqT}
 \sum_{T\in G_{M,i}}x^{{\rm xleaf}(T)}y^{{\rm yleaf}(T)}
 z^{{\rm zleaf}(T)}
 =  \sum_{T\in H_{M,i}}(xy)^{{\rm yleaf}(T)}(x+y)^{K+1-i-2{\rm yleaf}(T)} z^i,
 \end{equation}
 where $H_{M,i}$ denotes  the
 set of canonical  trees $T\in G_{M, i}$ without any unbalanced $y$-leaves.

For each $T\in H_{M,i}$, we define ${\rm Orbit}(T)$ to be the
set of Gessel trees $S$ that can be transformed to $T$
via the Foata-Strehl action.
   Clearly,  the number of $z$-leaves is invariant under any
   Foata-Strehl action.
   For any internal vertex $j$ of $T$ with an unbalanced $x$-leaf, there is a Gessel tree $S$
   such that $\psi_j(S)=T$. Let     ${\rm uxleaf}(T)$ and ${\rm bxleaf}(T)$ denote the
    numbers of unbalanced $x$-leaves and
    balanced $x$-leaves of $T$, respectively.
    Since there are no unbalanced $y$-leaves in $T$,
      \[ {\rm bxleaf}(T) =  {\rm yleaf}(T).\]
    Given that the total number of leaves in $T$ equals $K+1$,
     we find that    \[
                    2\, {\rm yleaf}(T)+ {\rm uxleaf}(T) ={\rm xleaf}(T) +{\rm yleaf}(T) = K+1-i,
          \]
          that is,
          \begin{equation}
               {\rm uxleaf}(T)  = K+1-i-2 {\rm yleaf}(T).
            \end{equation}
Therefore,
$$
 \sum_{S\in  {\rm Orbit}(T)}x^{{\rm xleaf(T)}}y^{{\rm yleaf(T)}}
z^{{\rm zleaf}(T)} = (xy)^{{\rm yleaf}(T)}(x+y)^{K+1-i-2{\rm yleaf}(T)} z^i.
$$
 Summing over all
 canonical Gessel trees in
 $H_{M,i}$ yields (\ref{eqT}), and hence the proof is complete.  \qed

\section{A context-free grammar approach}

 In this section, we present
 a context-free grammar approach to the partial
 $\gamma$-positivity of
 $C_{M}(x,y,z)$. Observe that the construction of
 Gessel trees implies a recursive
 formula for the computation of the Eulerian polynomials $C_M(x,y,z)$. For the case
 of
 Stirling permutations, the trivariate second order Eulerian
 polynomials $C_n(x,y,z)$ are defined by
 (\ref{C-n}). Dumont
\cite{Dumont-1980} deduced the
recursion
\begin{equation}
  C_n(x,y,z) = xyz
  \left(
 \frac{\partial}{\partial x}
 +
 \frac{\partial}{\partial y} +
 \frac{\partial}{\partial z}
  \right) C_{n-1}(x,y,z),
\end{equation}
where   $n\geq 1$ and $C_0(x,y,z)=x$,
see also Haglund-Visontai \cite{Haglund-Visontai-2012}.

\begin{theorem}
Let $n$ and $M$ be given as before. Set $M'=\{1^{k_1}, 2^{k_2}, \ldots, (n-1)^{k_{n-1}}\}$. Then
\begin{equation} \label{C-M-r}
  C_M(x,y,z) = xyz^{k_n-1}\left(
  \frac{\partial}{\partial x} +
  \frac{\partial}{\partial y} +
  \frac{\partial}{\partial z}\right) C_{M'}(x,y,z)
\end{equation}
with $C_\emptyset (x,y,z)=x$.
\end{theorem}

In the language of context-free grammars, for $k\geq 1$, define  the grammar
\begin{equation}
    G_k = \{ x \rightarrow xy z^{k-1}, \quad y \rightarrow xy
    z^{k-1}, \quad z \rightarrow xy
    z^{k-1} \}.
\end{equation}
Let $D_k$ denote the formal
derivative with respect to the
grammar $G_k$. Then the above
relation (\ref{C-M-r}) can be rewritten
as
\begin{equation}
   C_{M}(x,y,z) = D_{k_n}D_{k_{n-1}} \cdots
   D_{k_1}(x),
\end{equation}
where $D_{k_n}D_{k_{n-1}} \cdots
   D_{k_1}$ is meant to apply $D_{k_1}$ first,
   followed by the
   applications of $D_{k_2}$ and so on.
The grammatical expression
is informative to establish a connection to the
partial $\gamma$-positivity of $C_M(x,y,z)$.  Thanks to the
idea of the change of variables due to Ma-Ma-Yeh
\cite{Ma-Ma-Yeh-2019}, we set
\[ u=xy, \quad v=x+y \]
to get
\begin{equation}
    D_k(u) = D_k(xy) = D_k(x)y + x D_k(y) = xy(x+y) z^{k-1} = uv z^{k-1},
\end{equation}
\begin{equation}
    D_k(v) =D_k(x+y)=2xyz^{k-1}= 2u z^{k-1},
\end{equation}
and
\begin{equation}
    D_k(z) = xyz^{k-1}=u z^{k-1}.
\end{equation}
Now, for the variables $u,v, z$, the grammar $G_k$ can be recast as
\begin{equation}
\label{gk}
    G_k=\{ u \rightarrow uv z^{k-1},
     \quad v \rightarrow 2u z^{k-1},
     \quad z\rightarrow u z^{k-1}\}.
\end{equation}

The following theorem can be
viewed as an equivalent form of
Theorem \ref{mainthT} reconstructed
on a variation of Gessel trees, called
pruned Gessel trees,
along with a grammatical labeling
by using the variables $u$ and $v$.
It is worth mentioning that the grammatical labeling can be
thought as a guideline to generate the $\gamma$-coefficients. Indeed,
it drops a hint in search for
a Foata-Strehl action.

Intuitively, a pruned Gessel tree
is obtained from a Gessel tree by chopping off the
$x$-leaves and $y$-leaves. We have to say
that this viewpoint alone is of no substantial help.
Here comes the idea of characterizing pruned Gessel trees. Assume that
$T$ is a Gessel tree on $M$. Then
internal vertices of the
pruned Gessel tree obtained from $T$
are of the following four types.

 \noindent Type 1: $i$ has $k_i+1$ children, and neither
       the first nor the last child is a leaf.
       This means that $i$ has neither an $x$-leaf
        nor a $y$-leaf in $T$.

 \noindent Type 2: $i$ has $k_i$ children, and the first child
         of $i$ is not a leaf. This means that $i$
         has a $y$-leaf, but no $x$-leaf, in $T$.
         In this case, the vertex $i$ is associated with a label $y$.

\noindent Type 3:  $i$ has $k_i$ children, and the last child
         of $i$ is not a leaf. This means that $i$
         has an $x$-leaf, but no $y$-leaf, in $T$.
         In this case, the vertex is associated with a label
         $x$.

\noindent Type 4:  $i$ has $k_i-1$ children.  This means
       that $i$ has both an $x$-leaf and a $y$-leaf
        in $T$. In this case, the vertex $i$ is
        associated a label $xy$.

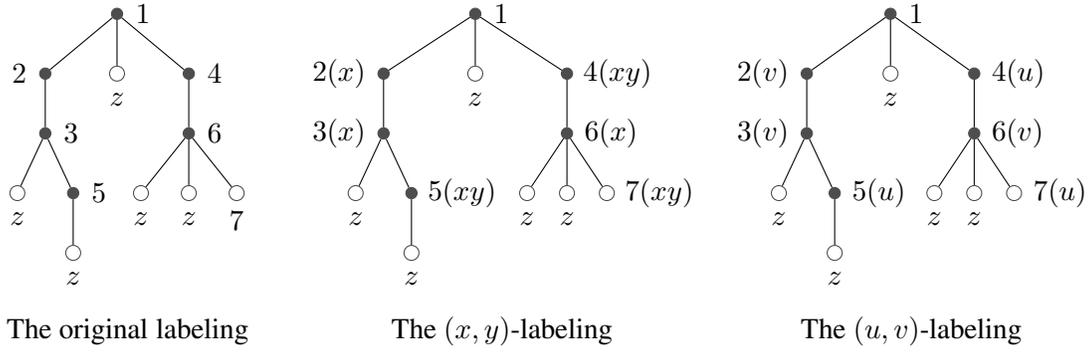
\begin{figure}
\begin{center}
\subcaptionbox*{\normalsize The original labeling}[.28\textwidth]{%
\begin{tikzpicture}[scale=0.53]
\node [tn1,label=0:{$1$}]{}[grow=down]
	[sibling distance=18mm,level distance=15mm]
        child {node [tn1,label=180:{$2$}](five){}
            [sibling distance=14mm,level distance=15mm]
            child{node [tn1,label=0:{$3$}]{}
            [sibling distance=14mm,level distance=15mm]
            child{node [tn,label=-90:{$z$}]{}}
                child{node [tn1,label=0:{$5$}](three){}
             [sibling distance=12mm,level distance=15mm]
            child{node [tn,label=-90:{$z$}](five2){}}
              edge from parent
                    node[right=2pt ]{}
                    }
                    }
            edge from parent
                node[left=2pt ]{}
            }
        child {node [tn,label=-90:{$z$}]{}
            edge from parent
                node[left=1pt ]{}
                }
        child {node [tn1,label=0:{$4$}]{}
            child{node [tn1,label=0:{$6$}](four){}
                [sibling distance=12mm,level distance=15mm]
            child{node [tn,label=-90:{$z$}]{}}
            child{node [tn,label=-90:{$z$}]{}}
            child{node [tn,label=-90:{$7$}]{}
              }
            edge from parent
                    node[right=2pt ]{}
            }
            edge from parent
                node[left=2pt ]{}
                };
\end{tikzpicture}
}
\subcaptionbox*{\normalsize The $(x,y)$-labeling}[.35\textwidth]{%
\begin{tikzpicture}[scale=0.53]
\node [tn1,label=0:{$1$}]{}[grow=down]
	[sibling distance=23mm,level distance=15mm]
        child {node [tn1,label=180:{$2(x)$}](five){}
            [sibling distance=14mm,level distance=15mm]
            child{node [tn1,label=180:{$3(x)$}]{}
            [sibling distance=14mm,level distance=15mm]
            child{node [tn,label=-90:{$z$}]{}}
                child{node [tn1,label=0:{$5(xy)$}](three){}
             [sibling distance=12mm,level distance=15mm]
            child{node [tn,label=-90:{$z$}](five2){}}
              edge from parent
                    node[right=2pt ]{}
                    }
                    }
            edge from parent
                node[left=2pt ]{}
            }
        child {node [tn,label=-90:{$z$}]{}
            edge from parent
                node[left=1pt ]{}
                }
        child {node [tn1,label=0:{$4(xy)$}]{}
            child{node [tn1,label=0:{$6(x)$}](four){}
                [sibling distance=10mm,level distance=15mm]
            child{node [tn,label=-90:{$z$}]{}}
            child{node [tn,label=-90:{$z$}]{}}
            child{node [tn,label=0:{$7(xy)$}]{}
              }
            edge from parent
                    node[right=2pt ]{}
            }
            edge from parent
                node[left=2pt ]{}
                };
\end{tikzpicture}
}
\subcaptionbox*{\normalsize The $(u,v)$-labeling}[.34\textwidth]{%
\begin{tikzpicture}[scale=0.53]
\node [tn1,label=0:{$1$}]{}[grow=down]
	[sibling distance=21mm,level distance=15mm]
        child {node [tn1,label=180:{$2(v)$}](five){}
            [sibling distance=14mm,level distance=15mm]
            child{node [tn1,label=180:{$3(v)$}]{}
            [sibling distance=14mm,level distance=15mm]
            child{node [tn,label=-90:{$z$}]{}}
                child{node [tn1,label=0:{$5(u)$}](three){}
             [sibling distance=12mm,level distance=15mm]
            child{node [tn,label=-90:{$z$}](five2){}}
              edge from parent
                    node[right=2pt ]{}
                    }
                    }
            edge from parent
                node[left=2pt ]{}
            }
        child {node [tn,label=-90:{$z$}]{}
            edge from parent
                node[left=1pt ]{}
                }
        child {node [tn1,label=0:{$4(u)$}]{}
            child{node [tn1,label=0:{$6(v)$}](four){}
                [sibling distance=10mm,level distance=15mm]
            child{node [tn,label=-90:{$z$}]{}}
            child{node [tn,label=-90:{$z$}]{}}
            child{node [tn,label=0:{$7(u)$}]{}
              }
            edge from parent
                    node[right=2pt ]{}
            }
            edge from parent
                node[left=2pt ]{}
                };
\end{tikzpicture}
}
\end{center}
\caption{A pruned canonical Gessel tree.}
\label{f5}
\end{figure}

Conversely, the four types of vertices
are sufficient for the generation of
pruned Gessel trees.
Figure \ref{f5} displays the
pruned tree of the canonical Gessel tree in Figure \ref{f4}
along with the $(x,y)$-labeling and the
$(u,v)$-labeling.

When restricted to pruned canonical Gessel trees,
Type 2 vertices are not allowed to show up.
This property enables us to
substitute the label $xy$ with $u$, and  the label $x$
with $v$.

The weight of
a pruned canonical
Gessel tree $T$, denoted by $w(T)$,
is defined to be the product of the
$(u,v)$-labels. As usual, the
empty product is meant to be
1.  For example, the weight of the
pruned Gessel tree in Figure \ref{f5}
equals $u^3v^3z^5$.
The $\gamma$-polynomial $\gamma_M(u,v,z)$, called the
$\gamma$-polynomial of
 $M$, is defined to be the
 generating function of the
 partial $\gamma$-coefficients
 in (\ref{gammaT}). To be more specific,
 let
 \begin{equation} \gamma_M(u,v,z) = \sum_{i=0}^{K-n}z^i \sum_{j=1}^{\lfloor
 (K+1-i)/2 \rfloor}
 \gamma_{M, i,j}u^jv^{K+1-i-2j}.
 \end{equation}

\begin{theorem}
Let $n$ and $M$ be given as before.  Then
 \begin{equation}
     \gamma_M(u,v,z) = \sum_{T\in P_M}  w(T),
 \end{equation}
 where   $P_M$  denotes the set of
pruned canonical Gessel trees on $M$.
\end{theorem}

For
the case of Stirling permutations,
the above grammar $G_k$ reduces to the grammar
given by Ma-Ma-Yeh \cite{Ma-Ma-Yeh-2019}, namely,
\begin{equation} \label{G-MMY}
 G= \{ u \rightarrow uv z,
     \quad v \rightarrow 2u z,
     \quad z\rightarrow u z\},
\end{equation}
where we have used $z$ in place of $w$ in \cite{Ma-Ma-Yeh-2019}.
The above grammatical labeling of
pruned canonical Gessel trees  for Stirling
permutations serve as an underlying combinatorial
structure for the grammar in (\ref{G-MMY}).

The following theorem asserts that
pruned canonical Gessel trees
on a multiset $M$ with
the $(u,v)$-labeling
can be generated in the same way as successively applying the
formal derivatives $D_{k_1}$,
$D_{k_2}$, $\ldots$, $D_{k_n}$ to $x$.

\begin{theorem}
Let $n$ and $M$ be given as before.
Then \begin{equation}\gamma_M(u,v,z)
     = D_{k_n}D_{k_{n-1}} \cdots
   D_{k_1}(x).
\end{equation}
\end{theorem}

The proof is in the same
lines as the argument in \cite{Chen-Fu-2022}  for
the grammatical labeling of 0-1-2
plane trees. Instead of presenting a  proof in full detail,
it suffices to focus on the
action corresponding to the
rule $v\rightarrow 2uz^{k-1}$
of the grammar $G_k$ as in
(\ref{gk}).

Assume that $T$ is a pruned
canonical Gessel tree on
$$
M'=\{1^{k_1}, 2^{k_2}, \ldots, (n-1)^{k_{n-1}}\}.
$$
As before, we have
$M=\{1^{k_1}, 2^{k_2},
\ldots, n^{k_{n}}\}$. Put $k=k_n$. Suppose that $T$
has a vertex $i$ with the label $v$. This means that
$i$ has $k$ children with the
last child not being a leaf.
There are two ways to generate
a pruned canonical Gessel tree
$S$ from $T$ by appending the
vertex $n$ to $T$ as a child of $i$.

\noindent Case 1. Make $n$ the first child of $i$.
Then $i$ will no longer has a label, and $n$ will
    be assigned the label $u$.
    Meanwhile, $n$ will have $k-1$ $z$-leaves.
    We see that this operation corresponds to
    the rule   $v\rightarrow uz^{k-1}$.

\noindent Case 2. Swap the
    first child and the last
    child along with their subtrees in the pruned
    canonical Gessel tree $S$
    obtained in Case 1. Then
    we get a pruned
    canonical Gessel tree on $M$. Summing up, these two
    cases correspond to the
    rule $v\rightarrow 2uz^{k-1}$.

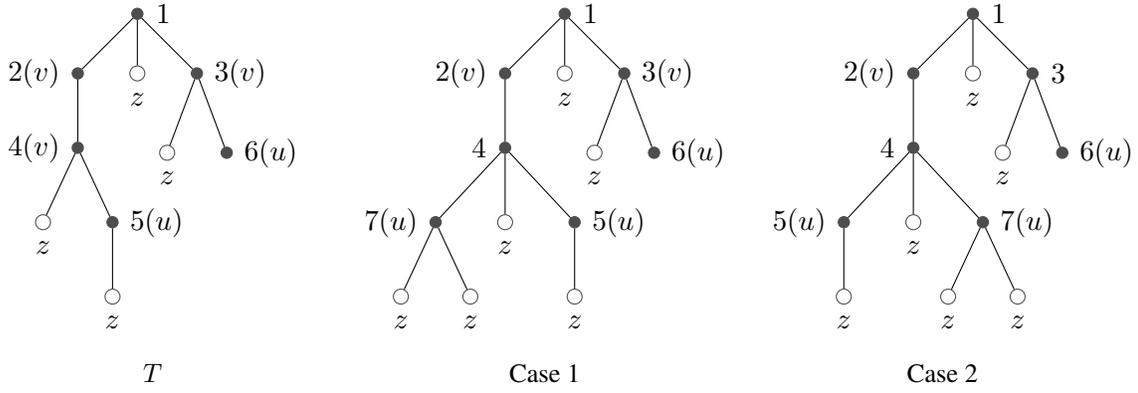
\begin{figure}
\begin{center}
\subcaptionbox*{$T$}[.3\textwidth]{%
\begin{tikzpicture}[scale=0.66]
\node [tn1,label=0:{$1$}]{}[grow=down]
	[sibling distance=12mm,level distance=12mm]
        child {node [tn1,label=180:{$2(v)$}](five){}
            [sibling distance=14mm,level distance=15mm]
            child{node [tn1,label=180:{$4(v)$}]{}
            [sibling distance=14mm,level distance=15mm]
            child{node [tn,label=-90:{$z$}]{}}
                child{node [tn1,label=0:{$5(u)$}](three){}
             [sibling distance=12mm,level distance=15mm]
            child{node [tn,label=-90:{$z$}](five2){}}
              edge from parent }
                    }
            edge from parent
                node[left=2pt ]{}
            }
        child {node [tn,label=-90:{$z$}]{}
            edge from parent
                node[left=1pt ]{}
                }
        child {node [tn1,label=0:{$3(v)$}]{}
           [sibling distance=12mm,level distance=16mm]
            child{node [tn,label=-90:{$z$}](three){}}
            child{node [tn1,label=0:{$6(u)$}](four){}
            edge from parent
                    node[right=2pt ]{}
            }
            edge from parent
                node[left=2pt ]{}
                };
\end{tikzpicture}
}
\subcaptionbox*{Case 1}[.36\textwidth]{%
\begin{tikzpicture}[scale=0.66]
\node [tn1,label=0:{$1$}]{}[grow=down]
	[sibling distance=12mm,level distance=12mm]
        child {node [tn1,label=180:{$2(v)$}](five){}
            [sibling distance=14mm,level distance=15mm]
            child{node [tn1,label=180:{$4$}]{}
            [sibling distance=14mm,level distance=15mm]
            child{node [tn1,label=180:{$7(u)$}]{}
            [sibling distance=14mm,level distance=15mm]
             child{node [tn,label=-90:{$z$}]{}}
             child{node [tn,label=-90:{$z$}]{}
             edge from parent
                node[left=2pt ]{}
                }
                edge from parent
                node[left=2pt ]{}
            }
            child{node [tn,label=-90:{$z$}]{}}
                child{node [tn1,label=0:{$5(u)$}](three){}
             [sibling distance=12mm,level distance=15mm]
            child{node [tn,label=-90:{$z$}](five2){}}
              edge from parent }
                    }
            edge from parent
                node[left=2pt ]{}
            }
        child {node [tn,label=-90:{$z$}]{}
            edge from parent
                node[left=1pt ]{}
                }
        child {node [tn1,label=0:{$3(v)$}]{}
           [sibling distance=12mm,level distance=16mm]
            child{node [tn,label=-90:{$z$}](three){}}
            child{node [tn1,label=0:{$6(u)$}](four){}
            edge from parent
                    node[right=2pt ]{}
            }
            edge from parent
                node[left=2pt ]{}
                };
\end{tikzpicture}
}
\subcaptionbox*{Case 2}[.31\textwidth]{%
\begin{tikzpicture}[scale=0.66]
\node [tn1,label=0:{$1$}]{}[grow=down]
	[sibling distance=12mm,level distance=12mm]
        child {node [tn1,label=180:{$2(v)$}](five){}
            [sibling distance=14mm,level distance=15mm]
            child{node [tn1,label=180:{$4$}]{}
            [sibling distance=14mm,level distance=15mm]
              child{node [tn1,label=180:{$5(u)$}](three){}
             [sibling distance=12mm,level distance=15mm]
            child{node [tn,label=-90:{$z$}](five2){}}
              edge from parent }
            child{node [tn,label=-90:{$z$}]{}}
                   child{node [tn1,label=0:{$7(u)$}]{}
            [sibling distance=14mm,level distance=15mm]
             child{node [tn,label=-90:{$z$}]{}}
             child{node [tn,label=-90:{$z$}]{}
             edge from parent
                node[left=2pt ]{}
                }
                edge from parent
                node[left=2pt ]{}
            }
                    }
            edge from parent
                node[left=2pt ]{}
            }
        child {node [tn,label=-90:{$z$}]{}
            edge from parent
                node[left=1pt ]{}
                }
        child {node [tn1,label=0:{$3$}]{}
           [sibling distance=12mm,level distance=16mm]
            child{node [tn,label=-90:{$z$}](three){}}
            child{node [tn1,label=0:{$6(u)$}](four){}
            edge from parent
                    node[right=2pt ]{}
            }
            edge from parent
                node[left=2pt ]{}
                };
\end{tikzpicture}
}
\end{center}
\caption{The two possibilities for the rule
$v\rightarrow 2uz^{k-1}$.
}
\label{f6}
\end{figure}

For example, for the pruned
canonical Gessel tree $T$ on
$M'=\{1^2, 2, 3^2, 4^2, 5^2, 6\}$
in Figure \ref{f6},
there are two ways to append the
vertex $7$ to $T$ as a child of $4$
to produce a pruned canonical
Gessel tree on $M=M' \cup \{ 7^3\}$.

It must be understood, however, that
the notion of pruned
canonical Gessel trees is
somewhat cosmetic, in the sense that whereas it
does not change anything in nature, it may have an effect on the impression.
In fact, we might
as well keep the original
$(x,y)$-leaves, and transport
the labels of the $(x,y)$-leaves to their parents. Nevertheless, the pruned version bears the advantage of
taking a simpler form
especially for permutations  where the structure
of 0-1-2 increasing plane trees comes into play.

To conclude this section, we claim that
the grammar $G_k$ captures all the possibilities of constructing pruned canonical
Gessel trees on $M$ from the ones on $M'$. It is only a matter of exercise to verify
that this is indeed the case.

\section{The partial $\gamma$-coefficients}

In this section, we demonstrate that
the combinatorial interpretation
of the partial $\gamma$-coefficients
of the Eulerian polynomials of Stirling permutations
on a multiset falls into the scheme of
canonical Gessel trees.

First, we need to translate the
defining property of an unbalanced $y$-leaf
into the language of Stirling permutations. This goal can be
achieved with the aid of Proposition
\ref{RGC} on the implications of
the $x$-leaves and $y$-leaves to
Stirling permutations.

Let $\sigma$ be a Stirling permutation of $M$.
Observe that a descent $i$ arises only when
$\sigma_i$ is the last occurrence. Similarly, an ascent
$j$ arises only when $\sigma_i$ is the first occurrence.
Assume that an index $i$ is a descent of $\sigma$, and assume that
$\sigma_p$ is the first occurrence of $\sigma_i$, where
$p\leq i$.
As an extension of the notion of a double descent of a
permutation of $[n]$, we say that $i$ is a double fall
of $\sigma$ if $i$ is a descent and $p-1$ is a descent
as well. Keep in mind the convention that
$\sigma_0=\sigma_{K+1}=0$.
Denote by   ${\rm dfall}(\sigma)$  the number of double falls  of $\sigma$.
For example, for the Stirling permutation  $\sigma=2533114664$,  the descents $2,9,10$  are not double falls,  whereas   the descent $4$ is  a double fall.

\begin{proposition}
\label{p-5-1}
Let $n$ and $M$ be given as before.
Assume that $\sigma$ is a
Stirling permutation on $M$ and $T$ is
the corresponding Gessel tree of $\sigma$. Then
an index $i$ is a double fall
of $\sigma$ if and only  if
the vertex $\sigma_i$ in $T$ has an unbalanced
$y$-leaf.
\end{proposition}

\noindent{\em Proof.}
Let $j=\sigma_i$, and let $S_j(\sigma)$
be the $j$-segment with the first occurrence
of $j$ at position $p$ and the last occurrence
of $j$ at position $q$. Moreover,
let
\begin{equation}\label{sj}
S_j(\sigma) = w_0 \, j \, w_1 \, j\,
\cdots \, w_{k_j -1} \, j \, w_{k_j}
\end{equation}
be the Gessel decomposition of $S_j(\sigma)$.
Assume that $j$ has an
unbalanced $y$-leaf. We proceed to show that $i$ is a
double fall of $\sigma$. The $y$-leaf indicates that
the last factor $w_{k_j}$ is empty, from which it
follows that $i$ is a descent, since $S_j(\sigma)$ is
surrounded by smaller elements at both ends. Suppose the first occurrence
of $j$ appears at position $p$. It remains to confirm that $p-1$
is also a descent. The $y$-leaf is
unbalanced means that the segment $w_0$
in (\ref{sj}) is nonempty. By the definition
of the Gessel decomposition, any element
in $w_0$ is greater than $j$, so that
$p-1$ must be a descent. This completes
the proof. \qed

Given the above characterization
of unbalanced $y$-leaves in terms
of Stirling permutations, we
obtain the following interpretation
of the partial $\gamma$-coefficients.

\begin{theorem}\label{mainthQ}
Let $n$, $M$ and $K$ be given as before.
   For $0\leq i\leq  K-n $ and $1\leq j\leq {\lfloor(K+1-i)/2 \rfloor}$, we have
 \begin{equation}\label{M-i-j}
     \gamma_{M, i,j}=|\{\sigma\in {Q}_{M}\mid {\rm plat}(\sigma)=i,\, {\rm des}(\sigma)=j,\, {\rm dfall}(\sigma)=0\}|.
     \end{equation}
\end{theorem}

Notice  that a double fall of a Stirling permutation of
a multiset boils down to a double descent of a permutation of $[n]$.
Therefore,  when specialized to $M=[n]$, Theorem \ref{mainthQ}
reduces to the $\gamma$-expansion (\ref{gammaA})
of the  bivariate Eulerian polynomials $A_n(x, y)$
due to  Foata and Sch\"uzenberger \cite{FS-70}.

\section{A theorem of Ma-Ma-Yeh}

In this section, we show
that our combinatorial interpretation of the partial
$\gamma$-coefficients given the preceding section reduces to a theorem of Ma-Ma-Yeh
\cite{Ma-Ma-Yeh-2019} subject to a restatement prompted by a symmetry consideration.

Recall that the partial
$\gamma$-coefficients $\gamma_{n,i,j}$ for $C_n(x,y,z)$ are defined by
\begin{equation}\label{c-n-gamma}
    C_n(x,y,z) = \sum_{i=1}^n z^i \sum_{j=0}^{\lfloor (2n+1-i)/2\rfloor}
      \gamma_{n, i, j} (xy)^j (x+y)^{2n+1-i-2j}.
\end{equation}

Let $\sigma\in  {Q}_{n}$. An index $1 \leq i \leq 2n$ is called
an {ascent-plateau}  if $\sigma_{i-1}< \sigma_{i}=\sigma_{i+1}$ and    a
{ descent-plateau} if $\sigma_{i-1}> \sigma_{i}=\sigma_{i+1}$.
Let ${\rm aplat}(\sigma)$ and ${\rm dplat}(\sigma)$
denote the numbers of ascent-plateaux and descent-plateaux, respectively.
The following finding is due to Ma-Ma-Yeh \cite{Ma-Ma-Yeh-2019}.

\begin{theorem} \label{thma}
For $n\geq 1$, $0\leq i\leq n$ and
$1\leq j\leq \lfloor (2n+1-i)/2\rfloor$,  we have
\begin{equation} \label{gnij}
\gamma_{n,i,j} = | \{ \sigma \in  {Q}_n \, | \,
  {\rm des} (\sigma) =i, \;
  {\rm aplat} (\sigma) = j, \;
  {\rm dplat} (\sigma)=0\} |.
  \end{equation}
\end{theorem}

For the case of Stirling permutations, that is, $M=[n]_2$, we
write $G_n$ for $G_M$.
By Theorem \ref{mainthT},
 $\gamma_{n, i, j}$ equals
 the number of trees $T\in G_n$
with $i$ $z$-leaves and $j$ $y$-leaves
but with no unbalanced $y$-leaves.
In the meantime, Theorem \ref{mainthQ} provides
a formulation resorting to the notion
of a double fall of a Stirling permutation.
In this situation, a double fall can be described as follows.
Let $n\geq 1$ and $\sigma\in Q_n$.
Assume that $i$ is a descent of $\sigma$. Let $\sigma_i=j$.
Then $\sigma_i$ must be the second occurrence of $j$ in
$\sigma$. Assume that $\sigma_p$ is the first
occurrence of $j$ in $\sigma$. Now, the descent
$i$ is called a double fall if $p-1$ is a descent of
$\sigma$ as well.

While there is no doubt that Theorem \ref{mainthQ} offers a
legitimate combinatorial statement, one could not
help wondering about the
connection to the
Ma-Ma-Yeh story. It turns out that the answer
lies in the symmetry of the polynomials $C_n(x,y,z)$.
As far as
Stirling permutations are concerned,  Theorem
\ref{mainthQ} can be reassembled with regards to
a slight twist of  canonical Gessel trees.
In fact, one only needs to interchange the roles of the $x$-leaves and the $y$-leaves.
Equivalently, we define a canonical Gessel tree on $[n]_2$,
or a canonical ternary increasing tree,
by imposing the constraint that there are no internal vertices having
 a $z$-leaf, but no $x$-leaf.

\begin{theorem}
For $n\geq 1$, the partial $\gamma$-coefficient
$\gamma_{n,i,j}$ for Stirling permutations as defined
by (\ref{c-n-gamma}) equals the number of canonical
ternary increasing trees on $[n]_2$ with $i$ $y$-leaves and $j$ $x$-leaves.
\end{theorem}

Notice that the above explanation of $\gamma_{n,i,j}$
can be directly justified in the same manner as the proof of Theorem
\ref{mainthT}.
The following proposition gives a characterization
of Stirling permutations corresponding to
canonical  ternary  increaing trees.

\begin{proposition} Let $n\geq 1$ and $\sigma \in Q_n$.
Let $T$ be the corresponding increasing ternary increasing tree
of $\sigma$.
Then $\sigma$ has a descent-plateau if and only if
$T$ contains an internal vertex having a $z$-leaf, but no $x$-leaf.
That is to say, $\sigma$ has no descent-plateaux
if and only if $T$ is canonical.
\end{proposition}

For example, the vertex $2$ in the ternary increasing tree
in Figure 1 has a $z$-leaf, but no $x$-leaf. The corresponding
Stirling permutation $\sigma=  33552217714664 $ has a  descent-plateau
$522$. On the other hand, Figure 7 furnishes a canonical
ternary increasing tree. The corresponding Stirling
permutation $\sigma= 22 33 55 1 77 1 4 66 4  $ contains no descent-plateaux.

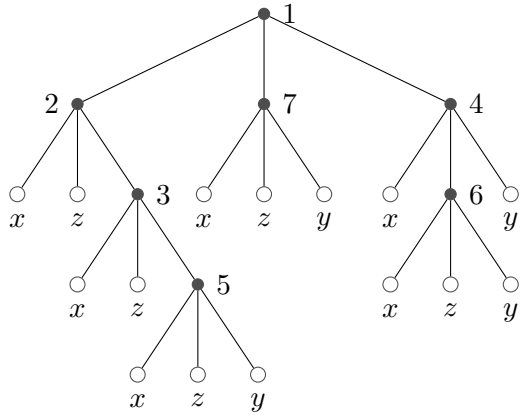
\begin{figure}[!ht]
\begin{center}
\begin{tikzpicture}[scale=0.8]
\node [tn1, label=0:{$1$}]{}[grow=down]
	[sibling distance=31mm,level distance=15mm]
        child {node [tn1,label=180:{$2$}](five){}
            [sibling distance=10mm,level distance=15mm]
            child{node [tn,label=-90:{$x$}](three){}}
            child{node [tn,label=-90:{$z$}](four){}}
             child{node [tn1,label=0:{$3$}]{}
            [sibling distance=10mm,level distance=15mm]
            child{node [tn,label=-90:{$x$}]{}}
            child{node [tn,label=-90:{$z$}](four){}}
         child{node [tn1,label=0:{$5$}](three){}
             [sibling distance=10mm,level distance=15mm]
                child{node [tn,label=-90:{$x$}](five1){}}
            child{node [tn,label=-90:{$z$}](five2){}}
            child{node [tn,label=-90:{$y$}]{}
              edge from parent
                    node[right=2pt ]{}
              }
              edge from parent
                    node[right=2pt ]{}
                    }
             edge from parent
                    node[right=2pt ]{}
                    }
            edge from parent
                node[left=2pt ]{}
            }
        child {node [tn1,label=0:{$7$}]{}
        [sibling distance=10mm,level distance=15mm]
             child{node [tn,label=-90:{$x$}]{}}
             child{node [tn,label=-90:{$z$}]{}}
             child{node [tn,label=-90:{$y$}]{}
                edge from parent
                node[left=1pt ]{}
                }
            edge from parent
                node[left=1pt ]{}
                }
        child {node [tn1,label=0:{$4$}]{}
           [sibling distance=10mm,level distance=15mm]
                child{node [tn,label=-90:{$x$}](three){}}
            child{node [tn1,label=0:{$6$}](four){}
                [sibling distance=10mm,level distance=15mm]
                child{node [tn,label=-90:{$x$}]{}}
            child{node [tn,label=-90:{$z$}]{}}
            child{node [tn,label=-90:{$y$}]{}
              edge from parent
                    node[right=2pt ]{}
              }
            edge from parent
                    node[right=2pt ]{}
            }
            child{node [tn,label=-90:{$y$}]{}
              edge from parent
                    node[right=2pt ]{}
              }
            edge from parent
                node[left=2pt ]{}
                }
        ;

\end{tikzpicture}

\end{center}
\caption{A canonical ternary
 increasing tree.
 }\label{f7}
\end{figure}

The proof of the above proposition is analogous to that of
Proposition \ref{p-5-1}. One more thing, we must add that
the statistic ${\rm aplat}(\sigma)$ reflects the number
of internal vertices in the corresponding
ternary increasing tree that have an $x$-leaf and a $z$-leaf
simultaneously. Thus we have reconfirmed the assertion of Ma-Ma-Yeh
in the context of canonical ternary increasing trees
in connection with Stirling permutations.

\vskip 5mm \noindent{\large\bf Acknowledgments.}
We wish to thank the referees for helpful suggestions.
This work was supported
by the National
Science Foundation of China.

\end{document}